%Paper: alg-geom/9404009
%From: Gregory Sankaran <G.K.Sankaran@pmms.cam.ac.uk>
%Date: Thu, 21 Apr 94 12:55 BST

% Plain TeX, version 3.1415
%
\def\CC {{\bf C}}

\def\HH {{\bf H}}

\def\QQ {{\bf Q}}

\def\ZZ {{\bf Z}}
\def\contin {\subseteq}
\def\Tilde{\widetilde}

\def\Bar{\overline}

\def\part#1#2 {{\partial {#1}/\partial {#2}}}

\def\Sp {\mathop{\rm Sp}\nolimits}

\def\Pic {\mathop{\rm Pic}\nolimits}

\def\im {\mathop{\rm Im}\nolimits}

\def\cE {{\cal E}}

\def\Sum{\sum\limits}

\def\tens{\otimes}

\def\pf{\noindent{\sl Proof: }}

\def\rationalmap{\mathrel{{\hbox{\kern2pt\vrule height2.45pt depth-2.15pt
 width2pt}\kern1pt {\vrule height2.45pt depth-2.15pt width2pt}
  \kern1pt{\vrule height2.45pt depth-2.15pt width1.7pt\kern-1.7pt}
   {\raise1.4pt\hbox{$\scriptscriptstyle\succ$}}\kern1pt}}}
\def\qed{\vrule width5pt height5pt depth0pt\par\smallskip}
\def\surj{\to\kern-8pt\to}
\outer\def\startsection#1\par{\vskip0pt
 plus.3\vsize\penalty-100\vskip0pt
  plus-.3\vsize\bigskip\vskip\parskip\message{#1}
   \leftline{\bf#1}\nobreak\smallskip\noindent}
\outer\def\camstartsection#1\par{\vskip0pt
 plus.3\vsize\penalty-150\vskip0pt
  plus-.3\vsize\bigskip\vskip\parskip
   \centerline{\it#1}\nobreak\smallskip\noindent}
\def\chain{\dot{\hbox{\kern0.3em}}}
\def\cochain{\d{\hbox{\kern0.3em}}}

\def\imic{\cong}
\outer\def\thm #1 #2\par{\medbreak
  \noindent{\bf Theorem~#1.\enspace}{\sl#2}\par
   \ifdim\lastskip<\medskipamount \removelastskip\penalty55\medskip\fi}
\outer\def\prop #1 #2\par{\medbreak
  \noindent{\bf Proposition~#1.\enspace}{\sl#2}\par
   \ifdim\lastskip<\medskipamount \removelastskip\penalty55\medskip\fi}
\outer\def\lemma #1 #2\par{\medbreak
  \noindent{\bf Lemma~#1.\enspace}{\sl#2}\par
   \ifdim\lastskip<\medskipamount \removelastskip\penalty55\medskip\fi}
\outer\def\corollary #1 #2\par{\medbreak
  \noindent{\bf Corollary~#1.\enspace}{\sl#2}\par
   \ifdim\lastskip<\medskipamount \removelastskip\penalty55\medskip\fi}
\def\deep #1 {_{\lower5pt\hbox{$#1$}}}

\raggedbottom
\baselineskip=24pt

\hoffset -.1in
\font\gothic=eufm10
\def \gothm {{\hbox{\gothic M}}}
\def \hf {{{1}\over{2}}}
\def \ap {{\cal A}_2(p)}
\def \BAp {{\bar{\cal A}_2(p)}}
\def \Ap {{\hat{\cal A}_2(p)}}
\def \Gp {\Gamma_{p^2}}
\def \TD#1 {{\Tilde{\Delta}_#1}}
\def \HD#1 {{\hat{\Delta}_#1}}
\def \HHD {{\skew{-8}\hat{\HD 0 }}}
\def \restr #1 {{\big\vert_{#1}}}
\def \phi {\varphi}
\centerline {\bf MODULI OF ABELIAN SURFACES WITH A $(1,p^2)$ POLARISATION}
\bigskip
\centerline{\bf V.A. Gritsenko \& G.K. Sankaran}
\bigskip
\noindent
The moduli space of abelian surfaces with a polarisation of type $(1,p^2)$ for
$p$ a prime was studied by O'Grady in [O'G], where it is shown that a
compactification of this moduli space is of general type if~$p\ge 17$. We
shall
show that in fact this is true if $p\ge 11$. Our methods overlap with those of
[O'G], but are in some important ways different. We borrow notation freely
from that paper when discussing the geometry of the moduli space.

\startsection 1. Methods

Let $\ap$ denote the moduli space of abelian
surfaces over $\CC$ with a polarisation of type $(1,p^2)$. We denote by
$\BAp$ the toroidal compactification of $\ap$ and by
$\Ap$ a partial desingularisation of $\BAp$ having only canonical
singularities.

To show that $\Ap$ is of general type (that is, roughly, that the
pluricanonical bundles have many sections), one chooses an
${\cal E}\in\Pic\Ap\tens\QQ$ such that $n\cal E$ (which is a bundle if $n$ is
sufficiently divisible) is not too far from the pluricanonical bundle
$nK_\Ap$, and such that the space of sections $H^0(n{\cal E})$ can be
calculated or at least estimated by some method. Then by knowing about the
geometry of $\Ap$ one can estimate the plurigenera, because the difference
between $nK_\Ap$ and $\cal E$ is known. In [O'G] the bundle used to play the
r\^ole of $n\cal E$ arises from pulling back powers of the Hodge bundle on the
moduli space $\Bar{\gothm}_2$ of semi-stable genus~$2$ curves via the map
$\Ap\to\Bar{\gothm}_2$ constructed there. Here, by contrast, we consider
$\ap$ as a Siegel modular variety, i.e., as a quotient of the Siegel
upper half-space by the paramodular group $\Gp$ (an arithmetic
subgroup of $\Sp (4,\QQ)$), and obtain a suitable bundle $n\cal E$ by
considering cusp forms of weight~$3n$ for~$\Gp$. A similar procedure
is adopted in [HS] for a different Siegel modular variety. But here we do not
use all cusp forms of weight~$3n$. Instead, we use the cusp form $f_2$
of weight~$2$ for $\Gp$, constructed by the first author in~[G], and
we consider modular forms of weight~$3n$ of the form
$f_nf_2^n$, where $f_n$ is a modular form of
weight~$n$. The bundle that results has fewer sections than the one arising
from all cusp forms but is much closer to~$nK_\Ap$ and this turns out to give
a better bound on the plurigenera for small~$p$. Hulek and the first author
have applied this idea to the situation of~[HS] and the improvement that
results is described in~[GH].

\startsection 2. Modular forms

If $t$ is a positive integer, the paramodular group is defined to be the
arithmetic subgroup of~$\Sp (4,\QQ)$
$$
\Gamma_t=\left\{\gamma\in\Sp (4,\QQ)\mid \gamma\in\pmatrix{
\phantom{t}\ZZ &\phantom{{{1}\over{t}}}\ZZ &\phantom{t}\ZZ &t\ZZ\cr
t\ZZ & \phantom{{{1}\over{t}}}\ZZ & t\ZZ & t\ZZ\cr
\phantom{t}\ZZ &\phantom{{{1}\over{t}}}\ZZ &\phantom{t}\ZZ &t\ZZ\cr
\phantom{t}\ZZ &{{1}\over{t}}\ZZ &\phantom{t}\ZZ &\phantom{t}\ZZ \cr
}\right\}.
$$
It acts on the Siegel upper half-plane
$$
\HH_2=\left\{Z=\pmatrix{\tau_1&\tau_2\cr
                       \tau_2&\tau_3\cr}
           \in M_{2\times 2}(\CC)\mid Z={}^TZ, \im Z >0\right\}
$$
by fractional linear transformations, i.e.,
$$
\pmatrix{A&B\cr C&D\cr}:Z\longmapsto (AZ+B)(CZ+D)^{-1}.
$$
The quotient $\Gamma_t\backslash \HH_2$ is a coarse moduli space for abelian
surfaces over $\CC$ with a polarisation of type~$(1,t)$. If $p$~is a prime and
$t=p^2$ we denote $\Gamma_t\backslash \HH_2$ by $\ap$, and let
$\BAp$ be the toroidal compactification and $\Ap$ the canonical partial
resolution described in [O'G]. If
we can show that $h^0(nK_\Ap)\sim n^3$ we shall have shown that $\Ap$
is of general type.

\prop 2.1 {\rm [G]} If $p\ge 11$ is a prime then there exists a nontrivial
cusp form $f_2$ of weight~$2$ for~$\Gp$.~\qed

It is not known whether such a cusp form exists for $\Gp$ if~$p\le 7$.

\prop 2.2 The space ${\cal M}_n^*(\Gp)$ of cusp forms of weight~$n$ for
$\Gp$ satisfies
$$
\dim {\cal M}_n^*(\Gp) = {{p^2(p^2+1)}\over{8640}}n^3+{\rm O}(n^2)
$$
for any prime~$p$.

\pf $\Gp$ is conjugate to a subgroup of $\Sp (4,\ZZ)$, namely
$$
\Gp'=\left\{\gamma\in\Sp (4,\ZZ)\mid \gamma\in\pmatrix{
\phantom{p}\ZZ &p\ZZ &\phantom{p}\ZZ &p\ZZ \cr
p\ZZ &\phantom{p}\ZZ &p\ZZ &\phantom{p}\ZZ \cr
\phantom{p}\ZZ &p\ZZ &\phantom{p}\ZZ &p\ZZ \cr
p\ZZ &\phantom{p}\ZZ &p\ZZ &\phantom{p}\ZZ \cr
}\right\}
$$
so we may as well work with $\Gamma'_{p^2}$. Exactly as in [HS] (cf. also [T])
we obtain, for~$l$ large
$$
\dim {\cal M}_n^*(\Gamma(l))\sim {{n^3}\over{8640}}[\Gamma(1):\Gamma(l)]
$$
and if $p|l$ then $\Gamma(l)\contin\Gp'$ and the cusp forms for $\Gp'$ are the
$\Gp'/\Gamma(l)$-invariant cusp forms for~$\Gamma(l)$. Using the Atiyah-Bott
fixed point theorem, as in~[T], we obtain
$$
\dim {\cal M}_n^*(\Gp')\sim {{2}\over{[\Gp':\Gamma(l)]}}
\dim {\cal M}_n^*(\Gamma(l))
$$
since there is a contribution from $\gamma=-I$ as well as from~$\gamma=I$. But
$$\eqalign{
{{2}\over{[\Gp':\Gamma(l)]}}\dim {\cal M}_n^*(\Gamma(l))
&\sim
{{2}\over{[\Gp':\Gamma(l)]}}\cdot{{n^3}\over{8640}}[\Gamma(1):\Gamma(l)]
\cr
&=
{{1}\over{[\Gp'/\pm I:\Gamma(l)]}}\cdot{{n^3}\over{8640}}[\Gamma(1):\Gamma(l)]
\cr
&=
{{n^3}\over{8640}}[\Gamma(1):\Gp'/\pm I]
\cr
&=
{{n^3}\over{8640}}\deg (\pi:\ap\to {\cal A}_2)
\cr
&=
{{p^2(p^2+1)}\over{8640}}n^3
\cr
}
$$
by [O'G, Lemma 2.1].~\qed

\startsection 3. Pluricanonical forms and extension to the boundary

Choose a cusp form $f_2$ of weight~$2$ for $\Gp$, $p\ge 11$; we can do
this in view of Proposition~2.1. Suppose $f_n$ is a modular form for
$\Gp$ of weight~$n$: then $\Phi=f_nf_2^n$ is a cusp form of
weight~$3n$. Let $\omega=d\tau_1\wedge d\tau_2\wedge d\tau_3$ be the standard
$3$-form on~$\HH_2$. The form $\Phi\omega^{\tens n}$ is invariant under
$\Gp$ and
therefore descends to give a pluricanonical form on $\ap$ except at
the branch locus of $\HH_2\to \ap$. If~$\Phi$ were a general element of
${\cal M}_{3n}^*(\Gp)$ we should expect this form to have logarithmic poles at
the boundary $\BAp\setminus \ap$, but because the
cusp form we have chosen is special these poles do not occur. That is because
$\Phi$~vanishes to high order (at least order~$n$) at the cusps.

\prop 3.1 The differential $3n$-form coming from $\Phi\omega^{\tens n}$
extends
over the generic point of each codimension~$1$ boundary component of
$\BAp$.

\pf According to [SC, Chapter~IV, Theorem~1] (see also [HS, Proposition~1.1]),
we need to check that in the Fourier-Jacobi expansion
$$
\Phi(Z)=\Sum_{m\ge 0}\theta^D_{m,\Phi}(u_D,t_D)\exp\{2\pi imz_D\}
$$
near the boundary component~$D$, the coefficients $\theta^D_{m,\Phi}$ vanish
for~$m<n$. But we can write the
expansion of $\Phi(Z)$ as a product of expansions
of $f_2(Z)$ and $f_n(Z)$: we have
$$
f_2(Z)=\Sum_{m>0}\theta^D_{m,f_2 }(u_D,t_D)\exp\{2\pi imz_D\}
$$
(with $\theta^D_{m,f_2 }(u_D,t_D)\equiv 0$ as $f_2$ is a cusp form), and
similarly for $f_n$. Hence
$$
\theta^D_{m,\Phi}=\Sum_{m_0+\cdots+m_n=m}\theta^D_{m_0,f_n }
        \prod_{i=1}^n\theta^D_{m_i,f_2 }
$$
which is zero if $m<n$ as then $m_i=0$ for some~$i\ge 1$.~\qed

$\BAp$ is smooth in codimension~$1$, but the quotient map
$\HH_2\to\ap$ is branched along two divisors (Humbert surfaces).
These are the divisors whose closures in $\BAp$ are denoted $\TD 1 $
and $\TD 2 $ in~[O'G]. At a point of $\HH_2$ lying over a general point of
$\TD 1 $ or $\TD 2 $ the isotropy group in $\Gp$ is $\ZZ/2$ and it acts by a
reflection, so $\BAp$ is smooth at a general point of $\TD 1 $ or~$\TD 2 $.

\corollary 3.2 If $n$ is sufficiently divisible then
$n(K_\BAp+\hf\TD 1 +\hf\TD 2 )$
is a bundle and if $\Phi=f_nf_2^n$ is a cusp form of the type described
above
then $\Phi\omega^{\tens n}$ determines an element of
$H^0\big(\BAp;n(K_\BAp+\hf\TD 1 +\hf\TD 2 )\big)$

\pf $\BAp$ has only quotient singularities, which are, in particular,
$\QQ$-Gorenstein. From the description of the action of $\Gp$ above a general
point of $\TD 1 $ or $\TD 2 $ it is clear that $\Phi\omega^{\tens n}$ acquires
poles
of order $\hf n$ along $\TD 1 $ and~$\TD 2 $.

\startsection 4. Obstructions from elliptic fixed points.

We take a canonical partial resolution $\phi:\Ap\to\BAp$, as
in~[O'G], adopting also the notations of [O'G, Definition 3.8] for (Weil)
divisors on~$\Ap$.

\prop 4.1 If $n$ is sufficiently divisible then $\Phi\omega^{\tens n}$
determines an element of
$$
H^0\big(\Ap;n(K_\Ap+\hf E'_1+\hf E''_1+\hf\HD 1 +\hf\HD 2
+(1-{{2}\over{p}})E_2-{{1}\over{4}}E'-{{1}\over{4}}E'')\big).
$$

\pf By Corollary 3.2., above, we have a
section of $\phi^*\bigl(n(K_\BAp+\hf\TD 1 +\hf\TD 2 )\bigr)$, and the formulae
for $\phi^*K_\BAp$, $\phi^*\TD 1 $ and $\phi^*\TD 2 $ given in [O'G]
provide the required expression.~\qed

We assume henceforth that $n$ is sufficiently divisible, so that everything we
have written so far is a bundle (in fact it is enough that $24p|n$). For the
rest of this section we assume that $p\ge 5$, as in [O'G], but we shall need
$p\ge 11$ in the end in order to apply Proposition~2.1.

Put $\cE=K_\Ap +\hf E'_1+\hf E''_1+\hf\HD 1 +\hf\HD 2 +(1-{{2}\over{p}})E_2
-{{1}\over{4}}E'-{{1}\over{4}}E''$. We want to make use of O'Grady's
calculations (and to avoid either resolving the singularities of $\Ap$ or
using adjunction and Riemann-Roch on singular varieties), so we express $\cE$
in terms of the pullback of the Hodge bundle $\lambda$ on $\Bar{\gothm}_2$
via the map $\pi\phi:\Ap\to\Bar{\gothm}_2$.

\lemma 4.2 In $\Pic \Ap \tens \QQ$ we have
$$
\cE = 3\phi^*\pi^*(\lambda) - {{1}\over{p}}\phi^*\pi^*(\Delta_1)
-{{p-1}\over{p}}\HD 0 -{{p-1}\over{p}}\HHD.
\eqno{(\star)\qquad}
$$

\pf See [O'G, Theorem 3.1].~\qed

We have
$$
\eqalign{
h^0(nK_\Ap)&\ge h^0(nK_\Ap-{{n}\over{4}}E'-{{n}\over{4}}E'')\cr
           &=  h^0\bigl(n\cE-{{n}\over{2}}E'_1-{{n}\over{2}}E''_1
   -{{n}\over{2}}\HD 1 -{{n}\over{2}}\HD 2 -n(1-{{2}\over{p}})E_2\bigr)
}
$$
so we want to estimate the five obstructions coming from $E'_1$, $E''_1$,
$\HD 1 $, $\HD 2 $ and $E_2$. Our $\cE$ plays the r\^ole of
$\phi^*\pi^*(\alpha_p\lambda)$ in [O'G]: by comparison, we have replaced
$\alpha_p=3-{{10}\over{p}}$ by~$3$, which makes some of the obstructing sheaves
more positive, but we also have vanishing to order
${{p-1}\over{p}}$ along the boundary components $\HD 0 $ and $\HHD$, which
makes them more negative. So much more negative, in fact, that we have the
following result.

\thm 4.3 All the obstructions vanish: that is, every section of $n\cE$ gives a
section of $nK_\Ap$ if $p\ge 5$ and $n$~is sufficiently divisible.

\pf We prove this in five steps, taking each obstruction separately. Only
$E'_1$ and $E''_1$ require much attention.

\noindent {\it 1. The obstruction from $E'_1$.}

We need (cf. [O'G], p.~146) to estimate
$h^0\bigl([n\cE-iE'_1]\restr {E'_1} \bigr)$ for $0\le i \le {{n}\over{2}}-1$.
Using ($\star$) we get
$$
\eqalign{
h^0\bigl([n\cE-iE'_1]\restr {E'_1} \bigr)
&= h^0\Bigl(3n\bigl[\phi^*\pi^*(\lambda)\bigr]\restr {E'_1}
-{{n}\over{p}}\bigl[\phi^*\pi^*(\Delta_1)\bigr]\restr {E'_1}
-n{{p-1}\over{p}}{\HD 0 }\restr {E'_1} -n{{p-1}\over{p}}{\HHD}\restr {E'_1}
-i{E'_1}\restr {E'_1} \Bigr)\cr
&\le h^0\Bigl(3n\bigl[\phi^*\pi^*(\lambda)\bigr]\restr {E'_1}
-n{{p-1}\over{p}}{\HD 0 }\restr {E'_1}
-i{E'_1}\restr {E'_1} \Bigr).\cr
}
$$
Replacing $\alpha_p$ by~$3$ in [O'G], Corollary 4.1 et seq., gives
$$
3n\bigl[\phi^*\pi^*(\lambda)\bigr]\restr {E'_1} -i{E'_1}\restr {E'_1}
 = 3i\Sigma + \bigl({{n}\over{4}}+4i\bigr)L - 3iG.
$$
As a set, $\HD 0 \cap E'_1$ consists of the fibres of
$\phi_2:E_1\to\Tilde\Gamma$ over the points where $\Tilde\Gamma$ meets the
boundary. These fibres are smooth, so $\HD 0  \restr {E'_1} $ is a positive
integer multiple of the general fibre~$L$. Hence
$$
\eqalign{
h^0\Bigl(3n\bigl[\phi^*\pi^*(\lambda)\bigr]\restr {E'_1}
-n{{p-1}\over{p}}{\HD 0 }\restr {E'_1}
-i{E'_1}\restr {E'_1} \Bigr)
&\le h^0\bigl(3i\Sigma
           + \bigl({{n}\over{4}}+4i-n{{p-1}\over{p}}\bigr)L - 3iG\bigr)\cr
&\le h^0\bigl(3i\Sigma
           + \bigl(({{1}\over{p}}-{{3}\over{4}})n+4i\bigr)L\bigr). \cr
}
$$
But
$$
\eqalign{
\Sigma\cdot\bigl(3i\Sigma+\bigl(({{1}\over{p}}-{{1}\over{4}})n+4i\bigr)L\bigr)
 &= ({{1}\over{p}}-{{3}\over{4}})n+i\cr
 &\le ({{1}\over{p}}-{{1}\over{4}})n-1\cr
}
$$
which is negative if $p\ge 5$, so there are no sections and no obstructions.

\noindent {\it 2. The obstruction from $E''_1$.}

The calculation here is very similar: this time it is $\HHD$ that plays a
part. We need to estimate
$h^0\bigl([n\cE-{{n}\over{2}}E'_1-iE''_1]\restr {E''_1} \bigr)$
for $0\le i \le {{n}\over{2}}-1$, and ($\star$) and [O'G, Corollary~4.2]
together give
$$
\eqalign{
h^0\bigl([n\cE-{{n}\over{2}}E'_1-iE''_1]\restr {E''_1} \bigr)
&\le h^0\Bigl(3n\bigl[\phi^*\pi^*(\lambda)\bigr]\restr {E''_1}
-{{n}\over{2}}E'_1\restr {E''_1}
-n{{p-1}\over{p}}{\HHD}\restr {E''_1}
-i{E''_1}\restr {E''_1} \Bigr)\cr
&\le h^0\Bigl(3i\Sigma + (4i-{{n}\over{4}})F
 - n{{p-1}\over{p}}{\HHD}\restr {E''_1} \Bigr)\cr
}
$$
(note that Lemma~4.7.(i) of [O'G] should read
$\phi^*\pi^*(\lambda)\restr {E''_1} \imic{{1}\over{12}}F$). As above,
${\HHD}\restr {E''_1} $ is a positive integer multiple of~$F$, so the
obstruction becomes
$h^0\bigl(3i\Sigma+\bigl(4i+({{1}\over{p}}-{{5}\over{4}})n\bigr)F\bigr)$. But
$$
\eqalign{
\Sigma\cdot\bigl(3i\Sigma+\bigl(4i+({{1}\over{p}}-{{5}\over{4}})n\bigr)F\bigr)
 &= i+({{1}\over{p}}-{{5}\over{4}})n\cr
 &\le ({{1}\over{p}}-{{3}\over{4}})n-1,\cr
}
$$
which is negative for all~$p$, so again there are no sections and no
obstructions.

\noindent {\it 3. The obstruction from $\HD 1 $.}

All we have to do is replace $\alpha_p$ by~$3$ in [O'G,~Corollary~4.4]. The
conclusion (Theorem~4.3 in [O'G]: the coefficient of $L'+L''$ should be
${{\alpha_p}\over{n}}+i$) that the obstruction vanishes is unaltered.

\noindent {\it 4. The obstruction from $\HD 2 $.}

Again changing $\alpha_p$ to~$3$ makes no difference: we simply get
$$
\displaylines{
\quad h^0\bigl(\bigl[n\cE-{{n}\over{2}}(E'_1+E''_1+\HD 1 )\bigr]\restr
{\HD 2 } \bigr)
    \hfill\cr
\qquad\qquad \le h^0\Bigl(3n\bigl[\phi^*\pi^*(\lambda)\bigr]\restr {\HD 2 }
-{{n}\over{2}}(E'_1+E''_1+\HD 1 )\restr {\HD 2 }
-n{{p-1}\over{p}}{\HD 0 +\HHD}\restr {\HD 2 }
-i{\HD 2 }\restr {\HD 2 } \Bigr)\hfill\cr
\qquad\qquad = 0 \qquad\hbox{ for all $i\ge 0$.}\hfill\cr
}
$$

\noindent {\it 5. The obstruction from $E_2$.}

This time the restriction of $\phi^*\pi^*(\lambda)$ is trivial, so $\alpha_p$
does not even appear in the calculation.~\qed

\thm 4.4 $\Ap$ is of general type if $p\ge 11$.

\pf This follows from Proposition~2.1, Proposition~2.2. amd Theorem~4.3.~\qed

In fact we have shown that
$h^0(nK_\Ap)\ge {{p^2(p^2+1)}\over{8640}}n^3 +{\rm O}(n^2)$ if~$p\ge 11$. We do
not really need the precise value of the leading coefficient (unless we really
want an asymptotic bound on the plurigenera, but O'Grady's bound is better
unless $p=11$ or $p=13$), because there are no obstructions to compare it with.
Instead, we have an explicit pluricanonical form.

\corollary 4.5 If $p\ge 11$ and $f_n$ is any modular form of weight~$n$ for
$\Gp$, then $f_n f_2^n \omega^{\tens n}$ gives an $n$-canonical
form on~$\Ap$.

\startsection References\par
\exhyphenpenalty100
\frenchspacing
\item{[SC]}A.~Ash, D.~Mumford, M.~Rapoport \& Y.~Tai, {\sl Smooth
Compactification of Locally Symmetric Varieties}, Math.\ Sci.\ Press,
Brookline, Mass., 1975.
\item{[G]}V.~Gritsenko, Irrationality of the moduli space of polarized
abelian surfaces, Bonn preprint MPI n-26/1994.
\item{[GH]}V.~Gritsenko \& K.~Hulek, Appendix to the paper
``Irrationality of the moduli space of polarized
abelian surfaces'',
Hannover preprint, 1994.
\item{[HS]}K.~Hulek \& G.~K. Sankaran, The Kodaira dimension of certain moduli
spaces of abelian surfaces,
Compositio\ Math.\ {\bf 90} (1994), 1--36.
\item{[O'G]}K.~O'Grady, On the Kodaira dimension of moduli spaces of
abelian surfaces, Compositio\ Math.\ {\bf 72} (1989), 121--163.
\item{[T]}Y.~Tai, On the Kodaira dimension of the moduli spaces of
abelian varieties, Inv.\ Math.\ {\bf 68} (1982), 425--439.

\raggedright
\vskip 1.0 true cm
V.~A. Gritsenko, Max-Planck-Institut f\"ur Mathematik,
Gottfried-Claren-Stra{\ss}e~26, D--53225~Bonn, Germany.
\vskip 1.0 true cm
G.~K. Sankaran, Department of Pure Mathematics and Mathematical Statistics,
16,~Mill~Lane, Cambridge CB2~1SB, England.

\end